\documentclass[11pt]{article}
\usepackage{tipa}
\usepackage{amsfonts}
\usepackage{amssymb}
\usepackage{amsmath}
\usepackage{amsbsy}
\usepackage{mathrsfs}
\usepackage{amssymb}
\usepackage{amsbsy}
\usepackage{amsfonts}
\usepackage{stmaryrd}
\def\proclaim#1{\par \bigskip\noindent {\bf #1}\bgroup\it\ }
\def\endproclaim{\egroup\par\bigskip}

\allowdisplaybreaks[4]

\newbox\TempBox \newbox\TempBoxA

\def\exp {\textsf{exp}}

\def\a.s{\textsf{a. s}}

\def\and{\textsf{and}}
\def\text#1{\mbox{\rm #1}}

\def\underwiggle 1{
\ifmmode\setbox\TempBox=\hbox{$ 1$}\else\setbox\TempBox=\hbox{ 1}\fi
\setbox\TempBoxA=\hbox to \wd\TempBox{\hss\char'176\hss}
\rlap{\copy\TempBox}\smash{\lower9pt\hbox{\copy\TempBoxA}} }

\begin{document}
\title{\bf  Weak Convergence to Stochastic Integrals Driven by  $\alpha-$Stable L\'{e}vy Processes \footnote{Project supported by the National Natural Science
Foundation of China (No. 10871177), and the Specialized Research
Fund for the Doctor Program of Higher Education (No.
20090101110020).}  }
\date{\today}
\author{Lin Zheng-Yan, Wang Han-Chao\footnote{Corresponding author, hcwang06@gmail.com.} \\
         Department of Mathematics, Zhejiang University, \\Hangzhou,
         China, 310027}

\maketitle $\mathbf{Abstract}$: We use the martingale convergence
method to get the weak convergence theorem on general functionals of
partial sums of independent heavy-tailed random variables. The
limiting process is the stochastic integral driven by
$\alpha-$stable L\'{e}vy process. Our method is very powerful to
obtain the limit behavior of heavy-tailed random variables.

 $\mathbf{Keyword}$: Weak convergence, martingale convergence, stochastic
 integral, $\alpha-$stable L\'{e}vy process, heavy-tailed.

\section{Introduction}
    $~~~~$Let $X_{n},n\ge 1,$ be independent and identically
    distributed (i.i.d.) random variables. When the distribution of
    $X_{1}$ is heavy-tailed, the limit behavior of stochastic
    processes which are related to $\{X_{n}\}$ are very important
    and interesting. In this paper, we will discuss the weak
    convergence  of following  processes:
          $$\sum_{i=2}^{[nt]}f(\sum_{j=1}^{i-1}(X_{n,j}-E(h(X_{n,j})))(X_{n,i}-E(h(X_{n,i}))),\eqno(1.1)$$
where $X_{n,j}=X_{j}/ b_{n}$ for some $b_{n}\rightarrow\infty$,
$f(x),h(x)$ are continuous functions.

This type of limit theorems is very important in probability theory,
mathematical statistics and econometrics, especially, it is a core
theory in the unit root model, which  is a hot topic in the
econometric theory (c.f. Phillips (1987 a,b), (2007)). In the unit
root theory, the limiting process of stochastic process sequence
like (1.1) is a  stochastic integral. In Ibragimov and Phillips
(2008), they studied the weak convergence of stochastic processes
like (1.1) when $X_{n}, n\ge1,$ are linear processes. Their theorems
are extension of unit root results. Lin and Wang (2010) studied the
same problems for causal processes.

In this paper, we extend these results to the heavy-tailed random
    variables. Heavy-tailed analysis is an interesting and important
    branch of probability, stochastic process and mathematical statistics.
Record-breaking insurance losses, financial log-returns,
transmission rates of files are examples of heavy-tailed phenomena.
     According to Rva\v{c}eva (1962), if the $X_{j},~j\ge 1,$ are i.i.d.,
     there exist $b_{n}>0$ and $c_{n}$ such that
      $$\frac{1}{b_{n}}\sum_{j=1}^{n}X_{j}-c_{n}\xrightarrow[]{d}\varsigma_{\alpha}\eqno(1.2)$$
for some non-degenerate $\alpha-$stable random variable
$\varsigma_{\alpha}$ with $\alpha\in (0,2)$ if and only if
 $X_{1}$ is regularly varying with index $\alpha\in (0,2)$.  After
then, a lot of authors studied the asymptotic behavior of
 independent or dependent heavy-tailed random variables.
 A detailed study of conditions for convergence of the partial sums of dependent stationary process to an infinite variance stable distribution was given in
 Bartkiewicz, Jakubowski, Mikosch and Wintenberger (2010), they
 also gave s survey for asymptotic distribution of partial sums of dependence heavy-tailed
 random variables.

 An extension of  (1.2) is its functional version
 for partial sum processes. In the other words, we consider
  the following
 processes:
           $$X_{n}(t)=\frac{1}{b_{n}}\sum_{i=1}^{[nt]-1}X_{i}-tc_{n},\eqno(1.3)$$
where $X_{n}(\cdot)$ is a random elements with values in the
Skorohod
 space $\mathbb{D}[0,1]$, i.e., the space of all function on $[0,1]$
 that are right-continuous and have left limits. Then the weak
 convergence of (1.3) is that for the probability
 measure on the space $\mathbb{D}[0,1]$. Many authors discussed this
  convergence. The point process method is a very powerful method to obtain this type
 of weak convergence. This method was  given in detailed by Resnick(1986). They combined the weak convergence of point processes with
 the
 continuous mapping theorem to obtain the results. They showed that
 that
the
 limiting process of (1.3) is $\alpha-$stable L\'{e}vy process if  $X_{1}$ is regularly varying with index $\alpha\in (0,2)$.
Davis and Hsing
 (1995) extends this result to dependent case.

 In this paper, we will discuss the weak convergence of stochastic processes
 (1,1). In fact, they  can  be seen as  the discretizations  of stochastic
 integrals. When $X_{1}$ is regularly varying with index $\alpha\in
 (0,2)$, we get that the limiting process of (1.1) is a stochastic
 integral driven by $\alpha-$stable L\'{e}vy process.

The weak convergence of (1.1) is
  interesting and difficulty from the theoretical point. If we use the point
process method to obtain the weak convergence, the summation
functional should be proved as a continuous functional respect to
the topology of Skorohod
 space $\mathbb{D}[0,1]$, and the limiting process should have a compound Poisson representation.
However,
  the summation functional  like (1.1) is difficult to be proved as a
  continuous functional in the Skorohod
 space $\mathbb{D}[0,1]$.  Moreover,  the stochastic
  integral driven by $\alpha-$stable L\'{e}vy process
don't have a compound Poisson representation.  The point process
method can not be
 used easily.

In this paper, we will use the stochastic calculus method to obtain
the result.  Since the limiting process is  a semimartingale, we
will use the predictable characteristics of semimartingale to
describe the asymptotic behavior of underlying processes. This is a
very common method in the study of classical stochastic analysis.
More details can be found in Jacod and Shiryaev
     (2003), which introduces
predicable characteristics to replace the three terms in the usual
case: the drift, variance of the Guassian part and the L\'{e}vy
measure, which characterize the distribution of the L\'{e}vy
process. By means of these three characteristics,  the tightness
criteria of semimartingale sequence is obtained. Furthermore, one
can identify the law of limiting process through the unique solution
of martingale problem related to the predicable characteristics. In
some special cases, the unique solution of a martingale problem can
be seen as a unique solution of stochastic differential equations
(for example, when the limiting process is a stochastic integral).
In this paper, we firstly compute the predicable characteristics of
stochastic integral,  and then we through the so-called martingale
convergence method to get the criteria conditions for weak
convergence. The
     martingale convergence method is also summarized  in Jacod and Shiryaev
     (2003). This method is based on the martingale characteristic of semimartingale. When the limiting process is a semimartingale, martingale
convergence method is  very powerful.

 The assumptions for obtaining main results are same as those in the point process method, our method is more simple
  than the point process method.   In Jacod and Shiryaev
     (2003), the authors use a same truncate function  to get the special semimartingale, and the predicable characteristics of      special martingale can determine the limit behavior of semimartingale.  However, we use a truncate function to get the special semimartingales for (1.1), and employ another
     different truncate function to deal with the limiting process.  We employ the core idea of the martingale convergence method to
    show the result.  Our method is a modification
  of the martingale convergence method. It is more convenience to
  verify the tightness conditions.
  As we know, the stochastic calculus method and martingale convergence method are not used in the
  asymptotic analysis of heavy-tailed phenomena in the previous study, our method may be  a new complement to the study
  of heavy-tailed phenomena. The similar
  procedure was  used in our previous paper, Lin and Wang (2010). Since
  the limiting process in that paper has no jumps, it is more simple
  than that in this paper. However, the jumps in the limiting process play a
  major role in the asymptotic analysis.

The remainder of this paper is organized as follows.
        Section 2 collects some basic tolls and notations to be used throughout this paper.
         In Section 3 and Section 4, we discuss the weak convergence to stochastic integrals driven by stable processes in the univariate and multivariate case respectively. Some
        discussion about the further research is given in Section 5.

\section{\bf{Preliminary}}
   \subsection{\bf{  Predictable Characteristics of Semimartingale and Convergence of Semimartingales.}}
$~~~~$ We follow the semimartingale theory as presented in Jacod and
     Shiryaev (2003). For our purpose, let
$\mathbb{R}_{+}=[0,+\infty)$ and $\mathbb{Z}=\{\cdots, -2, -1, 0, 1,
2,
      \cdots\}$. $(\Omega, \mathscr{F}, \mathbb{F}=(\mathscr{F}_{t})_{t\ge 1},
      P)$ is a filtered probability space. $X$ is a semimartingale
      defined on $(\Omega, \mathscr{F}, \mathbb{F}=(\mathscr{F}_{t})_{t\ge 1},
      P)$.  Set $h(x)$ is a continuous function satisfying $h(x)=x$ in a neighbourhood of $0$ and   $|h(x)|\le |x|1_{|x|\le
      1}$. Let
           $$
\left\{
\begin{array}{ll}
\check{X}(h)_{t}=\sum_{s\leq t}[\Delta X_{s}-h(\Delta X_{s})],\\
X(h)=X-\check{X}(h),
\end{array}
\right. $$  where $\Delta X_{s}=X_{s}-X_{s-}$. $X(h)$ is a special
semimartingale and we consider its canonical decomposition:
  $$X(h)=X_{0}+M(h)+B(h),\eqno(2.1)$$
where $M(h)$ is its local martingale part, $B(h)$ is its finite
variation part.

$\bf{Definition~1} $  (Jacod and
     Shiryaev (2003)) We call predictable characteristics of $X$
the triplet $(B,C,\nu)$ as follows:

  (1) $B$ is a predictable finite
variation process, namely the process $B=B(h)$ appearing in (2.1).

(2) $C=<M(h),M(h)>$ is a predictable process.

(3) $\nu$ is a predictable random measure on $\mathbb{R}_{+}\times
\mathbb{R}$, namely the compensator of the random measure $\mu^{X}$
associated to the jumps of $X$, $\mu^{X}$ is defined by
  $$\mu^{X}(\omega;dt,dx)=\sum_{s}1_{\{\triangle X_{s}(\omega)\neq 0\}}\varepsilon_{(s,\triangle X_{s}(\omega))}(dt,dx),\eqno(2.2)$$
where $\varepsilon_{a}$ denotes the Dirac measure at the point $a$,
which may be from different spaces.

$\bf{Definition~2} $  (Jacod and
     Shiryaev (2003)) Let $X$ be a c\`{a}dl\`{a}g
     process and let $\mathcal {H}$ be the $\sigma-$field generated
     by $X(0)$ and $\mathcal {L}_{0}$ be the distribution of $X(0)$.
     A solution to the martingale problem associated with $(\mathcal {H},
     X)$ and $(\mathcal {L}_{0},B,C,\nu)$ (denoted by $\varsigma (\sigma(X_{0}),X| \mathcal {L}_{0}, B, C,
      \nu)$) is a probability measure
     $P$ on $(\Omega, \mathscr{F})$ such that $X$ is a
     semimartingale on $(\Omega, \mathscr{F}, P)$ with predictable
     characteristics $(B,C,\nu)$.

The limit process $X=(X(s))_{s\ge 0}$ appearing in this paper is the
canonical process $X(s,\alpha)=\alpha(s)$ for the element
$\alpha=(\alpha(s))_{s\ge 0}$ of $D([0,1])$. In  other words, our
limit process is defined on the canonical space $(\mathbb{D}([0,1]),
\mathscr{D}([0,1]), \mathbf{D})$. For $a\ge 0$ and an element
$(\alpha(s), s\ge 0)$ of the Skorokhod space $\mathbb{D}([0,1])$,
define
            $$S^{a}(\alpha)=\inf(s:|\alpha(s)|\ge a~\text{or}~|\alpha(s-)|\ge a).$$

   In the paper, $\Rightarrow$ denotes weak convergence
   in an appropriate metric space, and $\xrightarrow[]{P}$ denotes
   convergence in probability.  $\mathbb{C}_{2}^{b}(R)$ denotes the set of all bounded continuous functions on $\mathbb{R}$ which are 0 around
   0. $\mathbb{C}_{1}^{b}(R)$ is a subclass of
   $\mathbb{C}_{2}^{b}(R)$ having only nonnegative functions, which
   contains all functions $g_{a}(x)=(a|x|-1)^{+}\wedge 1$ for all
   positive rationals $a$ and it is a convergence-determining class
   for the weak convergence induced by $\mathbb{C}_{2}^{b}(R)$.
For a finite variation process $A$,
   the total variation process of $A$ is denoted by $Var(A)$. For $K$ and $H$, $K\cdot H$ denotes the stochastic integral.
   The following propositions, provides the
   basis for the study of asymptotic properties of semimartingales, they can be found in Jacod and
     Shiryaev (2003).

     {\bf Proposition A} Let $\{X^{n},n\ge1\}$ be a sequence of c\`{a}dl\`{a}g processes, and suppose that for all $n,q\in N$, we have
     the  decomposition
             $$X^{n}=U^{nq}+V^{nq}+W^{nq}$$
satisfying that

 (i) the sequence $(U^{nq})_{n}$ is tight;

 (ii) the
sequence $(V^{nq})_{n}$ is tight and there is a sequence $(a_{q})$
of real numbers such that:
       $$\lim_{q\rightarrow\infty}a_{q}=0,~~\lim_{n\rightarrow\infty}P(\sup_{t\le 1}|\Delta V^{nq}|>a_{q})=0;$$

       (iii) for $\varepsilon>0$,
$$\lim_{q\rightarrow\infty}\limsup_{n\rightarrow\infty}P(\sup_{t\le 1}|W^{nq}|>\varepsilon)=0.$$

Then the sequence $(X^{n})$ is tight.

{\bf Proposition B} Let $Y^{n}$ be a c\`{a}dl\`{a}g process and
$M^{n}$ be a martingale on a same filtered probability space
$(\Omega, \mathscr{F}, \mathbb{F}=(\mathscr{F}_{t})_{t\ge 1},
      P)$.  Let $M$ be a
c\`{a}dl\`{a}g process defined on  the canonical space
$(\mathbb{D}([0,1]), \mathscr{D}([0,1]), \mathbf{D})$. Assume that

 (i) $(M^{n})$ is uniformly integrable;

 (ii) $Y^{n}\Rightarrow Y$ for some  $Y$ with law $\widetilde{P}=\mathscr{L}(Y);$

 (iii) $$M_{t}^{n}-M_{t}\circ (Y^{n})\xrightarrow[]{P},~~0\le t\le 1 $$
 Then the process $M\circ (Y)$ is a martingale under $\widetilde{P}$.

\subsection{Heavy-tailed  Random Variable and L\'{e}vy $\alpha-$Stable Process}

     $~~~~$In this subsection, we collect some facts, tools and notions
     about heavy-tailed random variables. Roughly speaking, a random
     variable $X$ heavy-tailed with index $\alpha\in (0,2)$ if there exists a positive
     parameter $\alpha$ such that
           $$P(X>x)\sim x^{-\alpha},~~x\rightarrow\infty.$$
Usually, people discuss a class of heavy-tailed random variables,
the so-called stable random variables which  will also be discussed
in this paper.

$\bf{Definition~3} $ A random variable $X$ is said to be
$\alpha$-stable if its characteristic function is given by
      $$E\exp{\{iuX\}}=\exp{\{iua_{\alpha}+\int(e^{iux}-1-iuh(x))\Pi_{\alpha}(dx)\}},$$
      $$a_{\alpha}=\begin{cases}
\beta\frac{\alpha}{1-\alpha},~~~~~\alpha\ne 1,\\
0,~~~~~~~~~~~\alpha=1,
\end{cases}$$
the index of stability $\alpha\in (0,2)$ and $\Pi_{\alpha}(dx)$ is
the L\'{e}vy measure.

In this paper, we use vague convergence to be assumption. Some
backgrounds on vague convergence  are given below. More details can
be found in Resnick (2007). More details can be found in Resnick
(2007). Let $E$ be a locally compact Hausdorff space with a
countable basis and $M_{p}(E)$ be the set of Radon measures on $E$
with values in $\mathbb{Z}_{+}$, where $\mathbb{Z}_{+}$ denotes the
set of positive integers. The space $M_{p}(E)$ is a Polish space
which is endowed with the topology of vague convergence. Recall that
for $\mu_{n},\mu\in M_{p}(E)$
             $$\mu_{n}\xrightarrow[]{v}\mu~~\text{iff}~~\mu_{n}(f)\rightarrow \mu(f)$$
for any $f\in C_{K}^{+}$, where $C_{K}^{+}$ is the class of
continuous functions with compact support. In this paper, we assume
$E=[-\infty,\infty]\backslash \{0\}$.

The stochastic integral, which will be considered, is driven by
     $\alpha-$stable L\'{e}vy process. It is a pure-jump process,
     in the other words, it can be presented as a point process. Let $X_{\alpha}(t)$ be a $\alpha-$stable L\'{e}vy
     Process,  $\Delta X_{\alpha}(t)= X_{\alpha}(t)-
     X_{\alpha}(t-)$, and
             $$\mu_{\alpha}(dt,dx):=\mu_{\alpha}(\omega, dt,dx)=\sum_{s}1_{\{\triangle X_{\alpha}(s)(\omega)\neq 0\}}\varepsilon_{(s,\triangle X_{\alpha}(s)(\omega))}(dt,dx).\eqno(2.3).$$
     We assume   the predictable compensator of
     $\mu_{\alpha}(dt,dx)$ is $ds\nu(dx)$, where $\nu(dx)$
is the L\'{e}vy measure of $X_{\alpha}(1)$. By the L\'{e}vy-It\^{o}
representation of L\'{e}vy process,
        $$X_{\alpha}(t)=\int_{0}^{t}\int h(x)(\mu_{\alpha}(ds,dx)-ds\nu(dx))+\int_{0}^{t}\int (x-h(x))\mu_{\alpha}(ds,dx),\eqno(2.4)$$
where $h(x)$ is a continuous truncate function.
\section{Convergence to Stochastic Integral Driven by a L\'{e}vy $\alpha-$Stable Process: The Univariate Case.}
      $~~~~$ In this section, we use the martingale convergence
      approach to obtain the weak convergence results for various
      general functionals of  partial sums of i.i.d. heavy-tailed
      random variables. The method of proof for these results is  new
      in the
      study of heavy-tailed analysis.

 Our main results are about the weak convergence of
   stochastic processes in $\mathbb{D}[0,1]$ with Skorohod $J_{1}$ topology.

 {\bf Theorem 1} Let $f:\mathbb{R}\rightarrow \mathbb{R}$ be a
 continuous differentiable function such that
          $$|f(x)-f(y)|\le K |x-y|^{a}\eqno(3.1)$$
 for some constants $K>0$, $a>0$ and all $x, y\in \mathbb{R}$.
 Suppose that $\{X_{n}\}_{n\ge 1}$ is a sequence of  i.i.d. random variables. Set
             $$X_{n,j}=\frac{X_{j}}{b_{n}}-E(h(\frac{X_{j}}{b_{n}}))\eqno(3.2)$$
for some $b_{n}\rightarrow\infty$. Define $\rho$ by
         $$\rho((x,+\infty])=px^{-\alpha},~~~~\rho([-\infty,-x))=qx^{-\alpha}\eqno(3.3)$$ for $x>0$, where $\alpha\in (0,1)$,
$0<p<1$ and $p+q=1$. Then
     $$(\sum_{i=1}^{[nt]}X_{n,i},\sum_{i=2}^{[nt]}f(\sum_{j=1}^{i-1}X_{n,j})X_{n,i})\Rightarrow (Z_{\alpha}(t),\int_{0}^{t}f(Z_{\alpha}(s-))dZ_{\alpha}(s)),\eqno(3.4)$$
in $\mathbb{D}[0,1]$, where $Z_{\alpha}(s)$ is an $\alpha-$stable
L\'{e}vy process with L\'{e}vy measure $\rho$ iff
      $$nP[\frac{X_{1}}{b_{n}}\in \cdot]\xrightarrow[]{v}\rho(\cdot) \eqno(3.5)$$
in $M_{p}(E)$.

{\bf Theorem 2} Let function $f$ be same as that in Theorem 1, and
$\{X_{n}\}_{n\ge 1}$ is a sequence of  i.i.d. random variables. Set
$$X^{\varepsilon}_{n,j}=\frac{X_{j}}{b_{n}}1_{\{|X_{j}|\ge
\varepsilon b_{n}\}}-E(h(\frac{X_{j}}{b_{n}}1_{\{|X_{j}|\ge
\varepsilon b_{n}\}}))$$ and
$$Z^{\varepsilon}_{\alpha}(t)=\int_{0}^{t}\int_{|x|>\varepsilon}
h(x)(\mu(ds,dx)-ds\nu(dx))+\int_{0}^{t}\int (x-h(x))\mu(ds,dx)$$ for
any $\varepsilon>0$ and  some $b_{n}\rightarrow\infty$. Define
$\rho$ as (3.3)
          for $\alpha\in [1,2)$,
Then
     $$(\sum_{i=1}^{[nt]}X^{\varepsilon}_{n,i}, \sum_{i=2}^{[nt]}f(\sum_{j=1}^{i-1}X^{\varepsilon}_{n,j})X^{\varepsilon}_{n,i})\Rightarrow
     (Z^{\varepsilon}_{\alpha}(t), \int_{0}^{t}f(Z^{\varepsilon}_{\alpha}(s-))dZ^{\varepsilon}_{\alpha}(s)), $$iff
      (3.5) stands.

{\bf Remark 1.}  Usually, $X_{1}$ is assumed to be mean zero and
symmetric random variable, but we don't have such assumptions.  It
can not be deal with easily through centralized and symmetric
procedure. When $\sum_{i=1}^{[nt]}X_{n,i}$ is integrator of
integral,
            $$\sum_{i=2}^{[nt]}f(\sum_{j=1}^{i-1}X_{n,j})(\sum_{k=1}^{i}X_{n,k}-\sum_{k=1}^{i-1}X_{n,k}),$$
            the mean part
 produces another stochastic processes through $f(\sum_{j=1}^{[nt]-1}X_{n,j})$,
            $$\sum_{i=2}^{[nt]}f(\sum_{j=1}^{i-1}X_{n,j})E(h(\frac{X_{j}}{b_{n}})).$$
            After centralized and symmetric
procedure, the limiting process of weak convergence maybe change.

{\bf Remark 2.} When $\alpha \in [1,2)$,
         $$\int_{0}^{1}x\rho(dx)=\infty.$$
it is different from the case of $\alpha \in (0,1)$ and  is more
difficult to obtain the same result as Theorem 1. We obtain a weaker
result.

Set
        $$Y_{n}(t)=\sum_{i=2}^{[nt]}f(\sum_{j=1}^{i-1}X_{n,j})X_{n,i},~~Y(t)=\int_{0}^{t}f(Z_{\alpha}(s-))dZ_{\alpha}(s),~~S_{n}(t)=\sum_{i=1}^{[nt]}X_{n,i}.$$
We want to prove
        $$H_{n}(t):=(Y_{n}(t),S_{n}(t))\Rightarrow H(t)=(Y(t),Z_{\alpha}(t)).$$

We firstly give some lemmas, which are the basis of the proof.

{\bf Lemma 1.} The predictable characteristics of
$(Y(t),Z_{\alpha}(t))$ are the terms $(B,C,\lambda)$ as follows:

 $$
\left\{
\begin{array}{ll}
B^{i}(t)=\left\{
\begin{array}{ll}\int_{0}^{t}\int
(h(f(Z_{\alpha}(s-)x)-f(Z_{\alpha}(s-)h(x))\nu(ds,dx),~~~~~~~~~~~~i=1,\\
0,~~~~~~~~~~~~~~~~~~~~~~~~~~~~~~~~~~~~~~~~~~~~~~~~~~~~~~~~~~~~~~~~~~~~~~~i=2,\end{array}
\right. \\
C^{ij}(t)=\left\{
\begin{array}{ll}\int_{0}^{t}\int
h^{2}(f(Z_{\alpha}(s-)x)\nu(ds,dx),~~~~~~~~~~~~~~~~~~~~~~~~~~i=1,j=1,\\
\int_{0}^{t}\int
h(f(Z_{\alpha}(s-)x)h(x)\nu_{n}(ds,dx),i=1,j=2,\text{or}~i=2,j=1,\\
\int_{0}^{t}\int h^{2}(x)
\nu(ds,dx),~~~~~~~~~~~~~~~~~~~~~~~~~~~~~~~~~~~~~~i=2,j=2,\end{array}
\right.\\
1_{G}*\lambda(ds,dx)=1_{G}(x,f(Z_{\alpha}(s-))x
)\nu(ds,dx)~\text{for}~\text{all}~G\in\mathbb{B}^{2},\end{array}
\right. $$ where $\nu(ds,dx)$ is  the compensator of the jump
measure of $Z_{\alpha}(t)$.

 {\bf Proof.} From (2.4) and $\nu(\{t\}\times dx)=0$, $Z_{\alpha}(t)$, $B(t)$
and $C^{22}(t)$ are obtained by Proposition 2.17 in Chapter 2 of
Jacod and Shiryaev(2003).

 Let
$\eta(ds,dx)$ be the jump random measure of $Y(t)$ and
$\lambda'(ds,dx)$ be the compensator of $\eta(ds,dx)$.

If $G$ is a Borel set in $R$, we have
        $$1_{G}*\lambda'(ds,dx)=1_{G}(f(Z_{\alpha}(s-))x)*\nu(ds,dx).$$
Set $z=f(Z_{\alpha}(s-))x$, then
\begin{eqnarray*}& &Y(t)-\int_{0}^{t}\int (z-h(z))\eta(ds,dz)\\
               &=&\int_{0}^{t}\int f(Z_{\alpha}(s-))h(x)(\mu(ds,dx)-\nu(ds,dx))+\int_{0}^{t}\int f(Z_{\alpha}(s-))(x-h(x))\mu(ds,dx)\\
               & &-\int_{0}^{t}\int (z-h(z))\eta(ds,dz)\\
               &=&\int_{0}^{t}\int f(Z_{\alpha}(s-))h(x)(\mu(ds,dx)-\nu(ds,dx))+\int_{0}^{t}\int f(Z_{\alpha}(s-))(x-h(x))\mu(ds,dx)\\
               & &-\int_{0}^{t}\int (f(Z_{\alpha}(s-))x-h(f(Z_{\alpha}(s-)x))\mu(ds,dx)\\
               &=&\int_{0}^{t}\int f(Z_{\alpha}(s-))h(x)(\mu(ds,dx)-\nu(ds,dx))\\ & &+\int_{0}^{t}\int
(h(f(Z_{\alpha}(s-)x)-f(Z_{\alpha}(s-)h(x))\mu(ds,dx)\\
               &=&\int_{0}^{t}\int f(Z_{\alpha}(s-))h(x)(\mu(ds,dx)-\nu(ds,dx))\\ & &+\int_{0}^{t}\int
h(f(Z_{\alpha}(s-))x)-f(Z_{\alpha}(s-))h(x)(\mu(ds,dx)-\nu(ds,dx))\\
      & & +\int_{0}^{t}\int
h(f(Z_{\alpha}(s-))x)-f(Z_{\alpha}(s-))h(x)\nu(ds,dx)\\
&=&\int_{0}^{t}\int
h(f(Z_{\alpha}(s-))x)(\mu(ds,dx)-\nu(ds,dx))\\
      & & +\int_{0}^{t}\int
h(f(Z_{\alpha}(s-))x)-f(Z_{\alpha}(s-))h(x)\nu(ds,dx),\end{eqnarray*}
which implies
         $$B_{t}^{1}=\int_{0}^{t}\int
h(f(Z_{\alpha}(s-))x)-f(Z_{\alpha}(s-))h(x)\nu(ds,dx),$$and the
martingale part of $Y_{t}$ is
     $$\int_{0}^{t}\int
h(f(Z_{\alpha}(s-))x)(\mu(ds,dx)-\nu(ds,dx)). \eqno (3.6)$$  Then we
can get $C^{11}$, $C^{12}$ and $C^{21}$. The lemma is proved.
$~~~~~~~~~~~~~~~~~~~~~~~~~~~~~~~~~~~~~~~~~~~~~~~~~~~~~~\blacksquare$

We set
  $$\mu_{n}(\omega;ds,dx)=\sum_{i=1}^{n}\varepsilon_{(\frac{i}{n},\frac{X_{i}(\omega)}{b_{n}})}(ds,dx),$$
then
$$\nu_{n}(\omega;ds,dx):=\sum_{i=1}^{n}\varepsilon_{(\frac{i}{n})}(ds)P(\frac{X_{i}}{b_{n}}\in
dx)$$ is the compensator of $\mu_{n}$ by the independent of
$\{X_{i}\}_{i\ge 1}$. Set
    $$\zeta_{n}(\omega;ds,dx)=\sum_{i=1}^{n}\varepsilon_{(\frac{i}{n},\frac{X_{i}(\omega)}{b_{n}}-c_{n})}(ds,dx),$$
we have
$$\varphi_{n}(\omega;ds,dx):=\sum_{i=1}^{n}\varepsilon_{(\frac{i}{n})}(ds)P(\frac{X_{i}}{b_{n}}-c_{n}\in
dx)$$ is the compensator of $\zeta_{n}(\omega;ds,dx)$, where
$c_{n}=E[h(\frac{X_{1}}{b_{n}})]$.

Firstly, we consider process $S_{n}(t)$. Introduce truncate function
$h_{n}(x)=h(x+c_{n})$.
            \begin{eqnarray*}S_{n}(t)&=& \sum_{i=1}^{[nt]}h(\frac{X_{i}}{b_{n}})+ \sum_{i=1}^{[nt]}(X_{n,i}-h(\frac{X_{i}}{b_{n}}))\\
                                     &=& \sum_{i=1}^{[nt]}(h(\frac{X_{i}}{b_{n}})-c_{n})+\sum_{i=1}^{[nt]}(\frac{X_{i}}{b_{n}}-h(\frac{X_{i}}{b_{n}}))\\
                                     &=&\int_{0}^{t}\int h(x)(\mu_{n}(ds,dx)-\nu_{n}(ds,dx))+\sum_{i=1}^{[nt]}(\frac{X_{i}}{b_{n}}-h(\frac{X_{i}}{b_{n}}))\\
                                     &=:& \widetilde{S}_{n}(t)+ \sum_{i=1}^{[nt]}(\frac{X_{i}}{b_{n}}-h(\frac{X_{i}}{b_{n}})).\end{eqnarray*}
The predictable characteristics of $\widetilde{S}_{n}(t)$ are
          $$B_{n}^{2}(t)=0,$$
          $$C_{n}^{22}(t)=\int_{0}^{t}\int h^{2}(x)\nu_{n}(ds,dx)-\sum_{s\le t}(\int h(x)\nu_{n}(\{s\},dx))^{2}.$$

For $Y_{n}(t)$, we have
            \begin{eqnarray*}Y_{n}(t)&=& \sum_{i=2}^{[nt]}h(f(\sum_{j=1}^{i-1}X_{n,j})\frac{X_{i}}{b_{n}})+ \sum_{i=2}^{[nt]}(f(\sum_{j=1}^{i-1}X_{n,j})X_{n,i}-h(f(\sum_{j=1}^{i-1}X_{n,j})\frac{X_{i}}{b_{n}}))\\
                                     &=& \sum_{i=2}^{[nt]}(h(f(\sum_{j=1}^{i-1}X_{n,j})\frac{X_{i}}{b_{n}})-E(h(f(\sum_{j=1}^{i-1}X_{n,j})\frac{X_{i}}{b_{n}})|\mathscr{F}_{i}))\\& &+\sum_{i=2}^{[nt]}(E(h(f(\sum_{j=1}^{i-1}X_{n,j})\frac{X_{i}}{b_{n}})|\mathscr{F}_{i})-f(\sum_{j=1}^{i-1}X_{n,j})E(h(\frac{X_{1}}{b_{n}})))\\
                                     & &+ \sum_{i=2}^{[nt]}(f(\sum_{j=1}^{i-1}X_{n,j})\frac{X_{i}}{b_{n}}-h(f(\sum_{j=1}^{i-1}X_{n,j})\frac{X_{i}}{b_{n}}))\\
                                     &=&\int_{0}^{t}\int h(f(\sum_{j=1}^{[ns]-1}X_{n,j})x)(\mu_{n}(ds,dx)-\nu_{n}(ds,dx))\\ & &+\int_{0}^{t}\int
                                     (h(f(\sum_{j=1}^{[ns]-1}X_{n,j})x)-f(\sum_{j=1}^{[ns]-1}X_{n,j})h(x))\nu_{n}(ds,dx)\\
                                     & &+\sum_{i=2}^{[nt]}(f(\sum_{j=1}^{i-1}X_{n,j})\frac{X_{i}}{b_{n}}-h(f(\sum_{j=1}^{i-1}X_{n,j})\frac{X_{i}}{b_{n}}))\\
                                     &=:& \widetilde{Y}_{n}(t)+\sum_{i=2}^{[nt]}(f(\sum_{j=1}^{i-1}X_{n,j})\frac{X_{i}}{b_{n}}-h(f(\sum_{j=1}^{i-1}X_{n,j})\frac{X_{i}}{b_{n}})).\end{eqnarray*}
The predictable characteristics of $\widetilde{Y}_{n}(t)$ are
          $$B_{n}^{1}(t)=\int_{0}^{t}\int
                                     (h(f(\sum_{j=1}^{[ns]-1}X_{n,j})x)-f(\sum_{j=1}^{[ns]-1}X_{n,j})h(x))\nu_{n}(ds,dx),$$
          $$C_{n}^{11}(t)=\int_{0}^{t}\int h^{2}(f(\sum_{j=1}^{[ns]-1}X_{n,j})x)\nu_{n}(ds,dx)-\sum_{s\le t}(\int h(f(\sum_{j=1}^{[ns]-1}X_{n,j})x)\nu_{n}(\{s\},dx))^{2},$$
          \begin{eqnarray*}C_{n}^{12}(t)&=&C_{n}^{21}(t)\\&=&
                                   \int_{0}^{t}\int h(f(\sum_{j=1}^{[ns]-1}X_{n,j})x)h(x)\nu_{n}(ds,dx)\\& &-\sum_{s\le t}(\int h(f(\sum_{j=1}^{[ns]-1}X_{n,j})x)\nu_{n}(\{s\},dx))(\int h(x)\nu_{n}(\{s\},dx)).\end{eqnarray*}
 {\bf Lemma 2}. Under (3.5),
  $$\int g(x)nF_{n}(dx)\rightarrow\int g(x)\rho(dx),~~~~~n\rightarrow\infty, \eqno(3.7)$$
for every continuous  $g\in \mathbb{C}_{2}^{b}(R)$,where
$F_{n}(x)=P(\frac{X_{1}}{b_{n}}\le x)$.

{\bf Proof.} From (3.5), we have
       $$\int h(x)nF_{n}(dx)\rightarrow\int h(x)\rho(dx),~~~~~n\rightarrow\infty,\eqno(3.8)$$
for every continuous compact support function $h$.

From (3.3), we can get that for any $\varepsilon>0$, there exists
$r>0$ such that $\rho((r,+\infty))+\rho((-\infty,-r))<\varepsilon$.

Set $B_{r}=[-r,r]$, we can find a continuous, compact support
function $g_{r}$, such that $1_{B_{r}}\le g_{r}\le 1$. Then
  \begin{eqnarray*}|\int g(x)nF_{n}(dx)-\int g(x)\rho(dx)|&\le&|\int g(x)nF_{n}(dx)-\int g(x)g_{r}(x)nF_{n}(dx)|+|\int g(x)g_{r}(x)nF_{n}(dx)\\& &-\int g(x)g_{r}(x)\rho(dx)|
  +|\int
g(x)g_{r}(x)\rho(dx)-\int g(x)\rho(dx)|\\
&\le& |\int g(x)g_{r}(x)nF_{n}(dx)-\int
g(x)g_{r}(x)\rho(dx)|+||g||(nF_{n}(B_{r}^{c})+\rho(B_{r}^{c})).\end{eqnarray*}
For $\varepsilon>0$, there exists $n_{0}$, such that as $n\ge
n_{0}$,
  $$|\int g(x)g_{r}(x)nF_{n}(dx)-\int
g(x)g_{r}(x)\rho(dx)|<\varepsilon.$$

From Theorem 3.2 (ii) in Resnick(2007), there exists $n_{1}$, as
$n\ge n_{1}$,
           $$|nF_{n}(B_{r}^{c})-\rho(B_{r}^{c})|< \varepsilon.$$
Then we have
  $$|\int g(x)nF_{n}(dx)-\int g(x)\rho(dx)|\le (3||g||+1)\varepsilon$$
as $n\ge \max\{n_{0},n_{1}\}$, which implies (3.7).
$~~~~~~~~~~~~~~~~~~~~~~~~~~~~~~~~~~~~~~~~~~~~~~~~~~~~~~~~~~~~~~~~~~~~~~~~~~~~~~~~~\blacksquare$

From (3.5), we can obtain
           $$\sum_{i=1}^{[nt]}X_{n,i}\Rightarrow Z_{\alpha}(t),\eqno(3.9)$$
by Corollary 7.1 in Resnick (2007).

So $\sum_{i=1}^{[nt]}X_{n,i}$ is related compact, in the other
words, $\sum_{i=1}^{[nt]}X_{n,i}$ is tightness.

By the tightness of $\sum_{i=1}^{[nt]}X_{n,i}$, we have that for any
$\varepsilon>0$, there are $n_{0}\in N$ and $K\in R^{+}$ with
  $$P(\sup_{t\le 1}|S_{n}(t)|>K)<\varepsilon~~\text{as}~~n\ge n_{0}.\eqno(3.10)$$

Since the convergence of $H_{n}(t)\Rightarrow H(t)$ is a local
property, it suffices to prove the Theorem 1 and 2 for
$f(S_{n}(t-))1_{[0,T]}$ and $f(Z_{\alpha}(t-))1_{[0,T]}$ for any
stopping time $T$.

We use $S_{n}^{C}$ and $S^{C}$ to replace $T$ in
$f(S_{n}(t-))1_{[0,T]}$ and $f(Z_{\alpha}(t-))1_{[0,T]}$
respectively, where $S_{n}^{C}=\inf(s:|S_{n}(s)|\ge
C~\text{or}~|S_{n}(s-)|\ge C)$.  As described in Pag\`{e}s (1986),
we can assume
        $$f(S_{n}(t-))\le C, f(Z_{\alpha}(t-))\le C\eqno(3.11)$$
identically for some constant $C$ in the following proof.

Let $\mathscr{K}$ be a compact subset of $R$ such that $|u|\le C$
for any $u\in \mathscr{K}$.

Set
$$1_{G}*\lambda_{n}(ds,dx)=1_{G}(x,f(\sum_{i=1}^{[ns]-1}X_{n,i})x
)\nu_{n}(ds,dx)~\text{for}~~G\in\mathbb{B}^{2}.$$
 {\bf Lemma 3} Under (3.5), we have that for $t>0$,
             $$Var[K*\lambda_{n}-(K*\lambda)\circ H_{n}]_{t}\xrightarrow[]{P}0,\eqno(3.12)$$
for every continuous  $K(x,u)\in \mathbb{C}_{2}^{b}(R\times
\mathscr{K})$ satisfying $K(x,u)=0$ for all $|x|\le \delta$, $u\in
\mathscr{K}$ for some $\delta>0$.

{\bf Proof.}  Since  this lemma is almost same as the Lemma IX 5.22
of Jacod and Shiryaev(2003). We verify thatthe assumptions of Lemma
IX 5.22 of Jacod and Shiryaev(2003) are satisfied.

At first, we show that for every continuous $g\in
 \mathbb{C}_{2}^{b}(R)$,
 $$Var[g*\nu_{n}-g*\nu]_{t}\rightarrow0~~~~\text{for}~t>0,\eqno(3.13)$$
which is   assumption (i) of the Lemma IX5.22 in Jacod and
Shiryaev(2003). In fact,
             $$\int_{0}^{t}\int g(x)\nu_{n}(ds,dx)=[nt]E(g(X_{n,1})),$$
and
            $$\int_{0}^{t}\int g(x)\nu(ds,dx)=t\int g(x)\rho(dx).$$
  We have
         $$Var[g*\nu_{n}-g*\nu]_{t}\le |\int g(x)nF_{n}(dx)-\int g(x)\rho(dx)|\frac{[nt]}{n}+|\frac{[nt]}{n}-t|\int g(x)\rho(dx),$$
  (3.13) is obtained by (3.7).

As proved in the Lemma IX5.22  in Jacod and Shiryaev(2003), we only
need prove (3.12) for $K(x,u)=g_{a}(x)g(x)R(u)$, where $R(u)$ is a
continuous function on $\mathscr{K}$, $g\in \mathbb{C}_{2}^{b}(R)$.

As described in the Lemma IX 5.22 of Jacod and Shiryaev(2003),
\begin{eqnarray*} & &Var[K*\lambda_{n}-(K*\lambda)\circ H_{n}]_{t}\\&\le& |R(f(\sum_{i=1}^{[nt]-1}X_{n,i}))|Var[gg_{a}*\nu_{n}-gg_{a}*\nu]_{t}+|R(f(S_{n}(t-)))-R(f(\sum_{i=1}^{[nt]-1}X_{n,i}))|\cdot (gg_{a}*\nu)_{t}\\
                                                                      &\le& ||R||Var[gg_{a}*\nu_{n}-gg_{a}*\nu]_{t}+||g|||R(f(S_{n}(t-)))-R(f(\sum_{i=1}^{[nt]-1}X_{n,i}))|\cdot (g_{a}*\nu)_{t}\end{eqnarray*}
We can get
                 $$||R||Var[gg_{a}*\nu_{n}-gg_{a}*\nu]_{t}\rightarrow 0$$
by (3.13). Since  $R(u)$ is a continuous function on $\mathscr{K}$,
$R(u)$ is uniformly continuous on  $\mathscr{K}$. For any
$\varepsilon>0$, there exists $\delta_{1}>0$, such that
$|y-y'|<\delta_{1}\Rightarrow |R(y)-R(y')|<\varepsilon$. Then we
have
\begin{eqnarray*}& &P(||g|||R(f(S_{n}(t-)))-R(f(\sum_{i=1}^{[nt]-1}X_{n,i}))|>\varepsilon)\\
&\le& P(|f(S_{n}(t-))-f(\sum_{i=1}^{[nt]-1}X_{n,i}))|>\frac{\delta_{1}}{||g||}) \\
                                                                                                              &\le & P(|X_{n,t}|>\frac{\delta_{1}}{||g||})\\
                                                                                                               &\le& 2\frac{\rho(\frac{\delta_{1}}{||g||}, \infty]}{n}\rightarrow 0\end{eqnarray*}
by the Lipschitz condition of $f$ and (3.5).

Then
              $$(||g||(R(f(S_{n}(t-)))-R(f(\sum_{i=1}^{[nt]-1}X_{n,i})))\xrightarrow[]{P}0. \eqno(3.14)$$
Since $g_{a}*\nu$ is a increase deterministic measure, and (3.14)
satisfies the assumption (ii) of the Lemma IX 5.22 in Jacod and
Shiryaev(2003),
             $$||g|||R(f(S_{n}(t-)))-R(f(\sum_{i=1}^{[nt]-1}X_{n,i}))|\cdot (g_{a}*\nu)_{t}\xrightarrow[]{P}0.$$  We complete the proof.
$~~~~~~~~~~~~~~~~~~~~~~~~~~~~~~~~~~~~~~~~~~~~~~~~~~~~~~~~~~~~~~~~~~~~~~~~~~~~~~~~~~~~~~~~~~~~~~~~~~~~~~~\blacksquare$

 {\bf Lemma 4} Under (3.5), we have
             $$Var[B_{n}^{1}-B^{1}\circ S_{n}]_{t}\xrightarrow[]{P}0~~~~\text{for}~t>0.\eqno(3.15)$$

{\bf Proof.} Let
           $$K(x,u)=h(ux)-uh(x).$$
We obtain the lemma by Lemma 3.
$~~~~~~~~~~~~~~~~~~~~~~~~~~~~~~~~~~~~~~~~~~~~~~~~~~~~~~~~~~~~~~~~~~~~~~~~~~~~~~~~~\blacksquare$

 {\bf Lemma 5} Under (3.5), we have
             $$Var[C_{n}^{ij}-C^{ij}\circ S_{n}]_{t}\xrightarrow[]{P}0~~~~\text{for}~t>0,\eqno(3.16)$$
where $i,j=1,2$.

 {\bf Proof.} We only prove the case of $i=1,j=1$. The other cases are similar.

  Although this lemma is different from Lemma 3, the method of proof is same as that of Lemma 3 through
 \begin{eqnarray*} &
&Var[h^{2}(f((\sum_{i=1}^{[nt]-1}X_{n,i})x)*\nu_{n}(ds,dx)-
h^{2}(f(Z_{\alpha}(s-)x)*\nu(ds,dx)\circ S_{n}]_{t}\\&\le& |f^{2}(\sum_{i=1}^{[nt]-1}X_{n,i})|Var[x^{2}*\nu_{n}-x^{2}*\nu]_{t}+|f^{2}(S_{n}(t-))-f^{2}(\sum_{i=1}^{[nt]-1}X_{n,i})|\cdot (x^{2}*\nu)_{t}\\
                                                                      &\le& C Var[x^{2}*\nu_{n}-x^{2}*\nu]_{t}+2C|f(S_{n}(t-))-f(\sum_{i=1}^{[nt]-1}X_{n,i})|\cdot (x^{2}*\nu)_{t}\end{eqnarray*}by $|h(x)|\le |x|1_{|x|\le1}$.

From (3.13) and (3.14),
         $$Var[h^{2}(f((\sum_{i=1}^{[nt]-1}X_{n,i})x)*\nu_{n}(ds,dx)-
h^{2}(f(Z_{\alpha}(s-)x)*\nu(ds,dx)\circ
S_{n}]_{t}\xrightarrow[]{P}0.$$

Hence in order to prove (3.16), It suffices to show

$$Var[\sum_{s\le
t}(\int
h(f(\sum_{j=1}^{[ns]-1}X_{n,j})x)\nu_{n}(\{s\},dx))^{2}]\xrightarrow[]{P}0\eqno(3.17)$$
which is equivalent to
        $$Var[\sum_{s\le t}(\int
h(f(\sum_{j=1}^{[ns]-1}X_{n,j})x)\nu_{n}(\{s\},dx))^{2}-\sum_{s\le
t}(\int h(f(Z_{\alpha}(s-))x)\nu(\{s\},dx))^{2}\circ
S_{n}]_{t}\xrightarrow[]{P}0,\eqno(3.18)$$ since $\nu(\{s\},dx))=0$.

However,
       $$Var[h(f((\sum_{i=1}^{[ns]-1}X_{n,i})x)*\nu_{n}(ds,dx)-
h(f(Z_{\alpha}(s-)x)*\nu(ds,dx)\circ S_{n}]_{t}\xrightarrow[]{P}0
\eqno(3.19)$$ can implies (3.18), and the proof of (3.19) is similar
to the above argument. We complete the
proof.$~~~~~~~~~~~~~~~~~~~~~~~~~~~~~~~~~~~~~~~~~~~~~~~~~~~~~~~~~~~~~~~~~~~~~~~~~~~~~~~~~~~~~~~~~~~~~~~~~~~~~~~~~~~~~~~~~~~~~~~~~~\blacksquare$

Set
$$1_{G}*\omega_{n}(ds,dx)=1_{G}(x,f(\sum_{i=1}^{[ns]-1}X_{n,i})x
)\varphi_{n}(ds,dx)~\text{for}~G\in\mathbb{B}^{2}.$$

{\bf Lemma 6} Under (3.5), we have that for $t>0$,
             $$Var[K*\omega_{n}-(K*\lambda)\circ S_{n}]_{t}\xrightarrow[]{P}0\eqno(3.20)$$
for every continuous  $K(x,u)\in \mathbb{C}_{2}^{b}(R\times
\mathscr{K})$ satisfying $K(x,u)=0$ for all $|x|\le \delta$, $u\in
\mathscr{K}$ for some $\delta>0$.

{\bf Proof.} Note that
           $$|c_{n}|\le E|\frac{X_{1}}{b_{n}}|1_{|X_{1}|\le b_{n}}=\int_{0}^{1}(P(|\frac{X_{1}}{b_{n}}|>y)-P(|\frac{X_{1}}{b_{n}}|>1))dy\rightarrow 0.$$
For $a\ne 0$,
         $$n(P(\frac{X_{i}}{b_{n}}-c_{n}<a)-P(\frac{X_{i}}{b_{n}}<a))\le nP(a-|c_{n}|\le \frac{X_{i}}{b_{n}}\le a+|c_{n}|)\rightarrow 0,$$
which implies
            $$nP[\frac{X_{1}}{b_{n}}-c_{n}\in \cdot]\xrightarrow[]{v}\rho(\cdot) \eqno(3.21)$$
by (3.5). From (3.21) and Lemma 2, we can get (3.20).
$~~~~~~~~~~~~~~~~~~~~~~~~~~~~~~~~~~~~~~~~~~~~~~~~~~~~~~~~~~~~~~~~~~~\blacksquare$

{\bf Remark 3} Based on the proof of Lemma 1-6, $\nu_{n}$ in Lemma 4
and 5 can be replaced by $\varphi_{n}$.

 {\bf Proof of Theorem 1 }Assume  (3.4) with $f(x)=x$  holds. From
Corollary 7.1 in Resnick (2007), we can get (3.5).

Assume that (3.5) holds.  we prove (3.4). The proof will be
presented in two steps.

 (a) We prove  the tightness of $H_{n}(t)$ by
using Theorem VI4.18 in Jacod and Shiryaev (2003).

The functions $\alpha\rightsquigarrow B_{t}(\alpha), C_{t}(\alpha),
      g*\lambda_{t}(\alpha)$ are Skorokhod-continuous on
      $\mathbb{D}(R)$ since the truncation function is continuous. Then $B_{n}(t), C_{n}(t),
      g*\omega_{n}(t)$ are C-tight by Lemmas 4-6.

From (3.5),
$$\mathscr{L}(S_{n}(t))\Rightarrow \mathscr{L}(Z_{\alpha}(t)).$$
It is means that $S_{n}(t)$ is tight. Note that
    $$\sum_{i=1}^{[nt]}\frac{X_{i}}{b_{n}}=S_{n}(t)+[nt]c_{n},$$
 and  $[nt]c_{n}\rightarrow\int_{0}^{t}\int h(x)\nu(ds,dx)$. Hence
$\sum_{i=1}^{[nt]}\frac{X_{i}}{b_{n}}$ is tight by Proposition A.

$$\lim
_{b\uparrow\infty}\lim\sup_{n}P(|x^{2}|1_{\{|x|>b\}}*\varphi
_{n}(t\wedge
         S_{n}^{a})>\varepsilon)=
         0\eqno(3.22)$$ for all $t>0$, $a>0$, $\varepsilon>0$
by the necessary part of Theorem VI4.18 in Jacod and Shiryaev
(2003).

 We have $$\lim
_{b\uparrow\infty}\lim\sup_{n}P(|x^{2}|1_{\{|x|>b\}}*\omega_{n}(t\wedge
         S_{n}^{a})>\varepsilon)=
         0$$  by (3.11), and hence
         $H_{n}(t)$ is tight.

(b) Identify  the limiting process. We need to prove that if a
subsequence, still denoted by
$\widetilde{P}^{n}=\mathscr{L}(H_{n})$, weakly  converges to a limit
$\widetilde{P}$ and  the semimartingale $H(t)$ has predicable
characteristics $(B,C,\lambda)$ under $\widetilde{P}$, we can
identify the limiting process. Since (3.1), the martingale problem
$\varsigma (\sigma(X_{0}),X| \mathcal {L}_{0}, B, C,
     \lambda)$ has unique solution by Theorem
6.13 in Applebaum (2009).

So our work is to prove the semimartingale $H$ has predicable
characteristics $(B,C,\lambda)$ under $\widetilde{P}$, in the other
words, to prove
      $$
h(f(Z_{\alpha}(s-))x)*(\mu(ds,dx)-\nu(ds,dx))\circ S_{n}(t), $$
      $$(h(f(Z_{\alpha}(s-))x)*(\mu(ds,dx)-\nu(ds,dx))\circ S_{n}(t))^{2}-C^{11}\circ S_{n}(t),$$
      $$g*\eta\circ S_{n}(t)-g*\lambda\circ S_{n}(t)~\text{for}~g\in C^{1}(R)$$
are local martingales under $\widetilde{P}$.

Since \begin{eqnarray*}& &\int_{0}^{t}\int
h(f(\sum_{j=1}^{[ns]-1}X_{n,j})x)(\zeta_{n}(ds,dx)-\varphi_{n}(ds,dx))-h(f(Z_{\alpha}(s-))x)*(\mu(ds,dx)-\nu(ds,dx))\circ
S_{n}(t)\\
&=&h(f((\sum_{i=1}^{[nt]-1}X_{n,i})x)*\varphi_{n}(ds,dx)-
h(f(Z_{\alpha}(s-)x)*\nu(ds,dx)\circ S_{n}(t),\end{eqnarray*}
 (3.19), Lemma 6 and Remark 3 implies
   $$\int_{0}^{t}\int
h(f(\sum_{j=1}^{[ns]-1}X_{n,j})x)(\zeta_{n}(ds,dx)-\varphi_{n}(ds,dx))-h(f(Z_{\alpha}(s-))x)*(\mu(ds,dx)-\nu(ds,dx))\circ
S_{n}(t)\xrightarrow[]{P}0. \eqno(3.23)$$

Set
    $$\widetilde{C}^{11}_{n}(t)=\int_{0}^{t}\int
h^{2}(f(\sum_{j=1}^{[ns]-1}X_{n,j})x)\varphi_{n}(ds,dx)-\sum_{s\le
t}(\int
h(f(\sum_{j=1}^{[ns]-1}X_{n,j})x)\varphi_{n}(\{s\},dx))^{2},$$ since
$$\mathscr{L}(S_{n}(t))\Rightarrow \widetilde{P},$$
 (3.11) implies that $C^{11}\circ S_{n}(t)\le C$, and Lemma 5, 6 implies that
 $P(\widetilde{C}^{11}_{n}(1)\ge
 C+1)\rightarrow 0$ as $n\rightarrow\infty$.
Set $T_{n}=\inf\{t: \widetilde{C}^{11}_{n}(t)>C+1\}$, we have
           $$\lim_{n\rightarrow\infty}P(T_{n}<1)=0.$$
so
  $$E(\sup_{0\le t\le 1}|\int_{0}^{t \wedge T_{n}}\int
h(f(\sum_{j=1}^{[ns]-1}X_{n,j})x)(\zeta_{n}(ds,dx)-\varphi_{n}(ds,dx))|^{2})\le
4E(\widetilde{C}_{n}^{11}(T_{n}))\eqno(3.24)$$ by Doob's inequality.

Since
$$\int_{0}^{t}\int
h(f(\sum_{j=1}^{[ns]-1}X_{n,j})x)(\zeta_{n}(ds,dx)-\varphi_{n}(ds,dx))$$
is a local martingale,  (3.23) and (3.24) imply that
$$h(f(Z_{\alpha}(s-))x)*(\mu(ds,dx)-\nu(ds,dx))\circ
S_{n}(t)$$ is a local martingale under $\widetilde{P}$ by
Proposition B.

A simple computation obtain that
\begin{eqnarray*}&
&(\int_{0}^{t}\int
h(f(\sum_{j=1}^{[ns]-1}X_{n,j})x)(\zeta_{n}(ds,dx)-\varphi_{n}(ds,dx)))^{2}\\&
& -(h(f(Z_{\alpha}(s-))x)*(\mu(ds,dx)-\nu(ds,dx))\circ
S_{n}(t))^{2}+ C^{11}\circ S_{n}(t)-\widetilde{C}_{n}^{11}(t)\\
&=&(\int_{0}^{t}\int
h(f(\sum_{j=1}^{[ns]-1}X_{n,j})x)(\zeta_{n}(ds,dx)-\varphi_{n}(ds,dx))\\&
&\cdot(\int_{0}^{t}\int
h(f(\sum_{j=1}^{[ns]-1}X_{n,j})x)(\zeta_{n}(ds,dx)-\varphi_{n}(ds,dx))-h(f(Z_{\alpha}(s-))x)*(\mu(ds,dx)-\nu(ds,dx))\circ
S_{n}(t)))\\
& &+ (h(f(Z_{\alpha}(s-))x)*(\mu(ds,dx)-\nu(ds,dx))\circ
S_{n}(t)))\\
& &\cdot(\int_{0}^{t}\int
h(f(\sum_{j=1}^{[ns]-1}X_{n,j})x)(\zeta_{n}(ds,dx)-\varphi_{n}(ds,dx))-h(f(Z_{\alpha}(s-))x)*(\mu(ds,dx)-\nu(ds,dx))\circ
S_{n}(t))\\
& &+C^{11}\circ S_{n}(t)-\widetilde{C}_{n}^{11}(t).\end{eqnarray*}
 We have that
$$\int_{0}^{t\wedge T_{n}}\int
h(f(\sum_{j=1}^{[ns]-1}X_{n,j})x)(\zeta_{n}(ds,dx)-\varphi_{n}(ds,dx))$$
  is uniformly integrable by (3.24),
  thus
  \begin{eqnarray*}&
&(\int_{0}^{t}\int
h(f(\sum_{j=1}^{[ns]-1}X_{n,j})x)(\zeta_{n}(ds,dx)-\varphi_{n}(ds,dx)))^{2}\\&
& -(h(f(Z_{\alpha}(s-))x)*(\mu(ds,dx)-\nu(ds,dx))\circ
S_{n}(t))^{2}+ C^{11}\circ
S_{n}(t)-\widetilde{C}_{n}^{11}(t)\xrightarrow[]{P}
0~~~~~~~~~~~(3.25)\end{eqnarray*} by (3.23), Lemma 5 and Remark 3.

By Lemma VII 3.34 of Jacod and Shiryaev (2003),
            $$E(\sup_{0\le t\le 1}|\int_{0}^{t \wedge T_{n}}\int
h(f(\sum_{j=1}^{[ns]-1}X_{n,j})x)(\zeta_{n}(ds,dx)-\varphi_{n}(ds,dx))|^{4})\le
K_{1}[E(\widetilde{C}_{n}^{11}(T_{n}))^{2}]^{\frac{1}{2}}+K_{2}E(\widetilde{C}_{n}^{11}(T_{n}))^{2}\eqno(3.26)$$
where $K_{1}$ and $K_{2}$ are constants.

Since
$$(\int_{0}^{t}\int
h(f(\sum_{j=1}^{[ns]-1}X_{n,j})x)(\zeta_{n}(ds,dx)-\varphi_{n}(ds,dx)))^{2}-\widetilde{C}_{n}^{11}(t)$$
is a local martingale,  (3.25) and (3.26) implies
$$(h(f(Z_{\alpha}(s-))x)*(\mu(ds,dx)-\nu(ds,dx))\circ
S_{n}(t))^{2}-C^{11}\circ S_{n}(t)$$ is local martingale under
$\widetilde{P}$ by Proposition B.

For $$g*\eta\circ S_{n}(t)-g*\lambda\circ S_{n}(t)~\text{for}~g\in
C^{1}(R),$$ we can get the similar conclusion by Lemma 6. We
complete the proof.
$~~~~~~~~~~~~~~~~~~~~~~~~~~~~~~~~~~~~~~~~~~~~~~~~~~~~~~~~~~~~~~~~~~~\blacksquare$

The proof of Theorem 2 is similar  except minor changes, we omit it
here.
\section{Convergence to Stochastic Integral Driven by L\'{e}vy $\alpha-$Stable Process: The Multivariate Case.}
$~~~~$ In this section, we use the similar method  to obtain the
weak convergence  for various
      general functionals of partial sums of i.i.d. heavy-tailed
      random vectors. Since the idea and method is similar, the proof
      are not given.

 {\bf Theorem 3} Let $f:\mathbb{R}\rightarrow \mathbb{R}$ be a
 continuous differentiable function such that
          $$|f(x)-f(y)|\le K |x-y|^{a}$$
 for some constants $K>0$, $a>0$ and all $x,y\in \mathbb{R}$.
 Suppose that $\{\xi_{n}\}_{n\ge 1}:=\{(\xi_{n}^{1},\xi_{n}^{2})\}_{n\ge 1}$ are i.i.d. random vectors. Set
             $$\xi_{n,j}=\frac{\xi_{j}}{b_{n}}-E(h(\frac{\xi_{j}}{b_{n}}))$$
for some $b_{n}\rightarrow\infty$.  Then
     $$(\sum_{i=1}^{[nt]}\xi_{n,i},\sum_{i=2}^{[nt]}f(\sum_{j=1}^{i-1}\xi^{1}_{n,j})\xi^{2}_{n,i})\Rightarrow (Z_{\alpha}(t),\int_{0}^{t}f(Z_{\alpha}^{1}(s-))dZ_{\alpha}^{2}(s))$$
in $\mathbb{D}[0,1]$, where $Z_{\alpha}(s)$ is a 2-dimesional
$\alpha-$stable L\'{e}vy Process with L\'{e}vy measure $\nu$ iff
$\xi_{1}$ is a  random vector satisfying the usual multivariate
regular variation condition with exponent $\alpha$ and
      $$nP[\frac{\xi_{1}}{b_{n}}\in \cdot]\xrightarrow[]{v}\nu(\cdot) \eqno(4.1)$$
 in $M_{p}(E_{2})$, where $E_{2}=[-\infty,\infty]\setminus \{0\}\otimes [-\infty,\infty]\setminus \{0\}.$

\section{Discussion.}

$~~~~$In this paper, we use a continuous function $h(x)$ for
technical convenience. In fact, we can take $h(x)=x1_{|x|\le 1}$ to
replace the continuous function through small change.

We only discuss  independent random variables. It  will be more
complex for dependence case. Recently, a lot of authors discussed
the functional limit theorems for
       $$X_{n}(t)=\frac{1}{b_{n}}\sum_{i=1}^{[nt]-1}X_{i}-tc_{n}$$
under dependence assumption (see Balan and Louhichi (2009),
Tyran-Kkami\'{n}ska(2010 a,b)). They employed the point process
method to deal with dependence. We hope that the method used in this
paper will be useful for study of the dependence heavy-tail random
variables.

\section*{Reference}

 Applebaum D. (2009). {\em L\'{e}vy Processes and Stochastic Calculus. 2nd edition}. Cambridge Press.\\
Balan R, Louhichi S. (2009). Convergence of point processes with
weakly dependent points. {\em Journal of Theoretical Probability
22}, 955-982.\\
Bartkiewicz K, Jakubowski A, Mikosch T, Wintenberger O. (2010).
Stable limits for sums of dependent infinite variance random
variables. Forthcoming in {\em Probability Theory and Related Fields}. \\
Davis R.A,  Hsing T. (1995). Point process and partial sum
convergence for weakly dependent random variables with infinite
variance. {\em Annals of Probability 23}, 879-917.\\
 Ibragimov
R, Phillips P. (2008). Regression asymptotics using martingale
convergence methods. {\em Econometric Theory 24},
888-947.\\
 Jacod J, Shiryaev AN.
(2003). {\em Limit Theorems for Stochastic
   Processes. } Springer.\\
Lin Z-Y, Wang H-C. (2010). On convergence to stochastic integrals.
 Arxiv preprint arXiv:1006.4693, 2010 .\\
Pag\`{e}s G. (1986) Un th\'{e}or\`{e}mes de convergence fonctionnel pour les int\'{e}grals stochastiques. {\em S\'{e}minaire de Proba. XX. Lecture Notes in Mathematics 1204, 572-611}.\\
Phillips P.C.B. (1987 a). Time-series regression with a unit root.
{\em Econometrica 55}, 277-301.\\
Phillips P.C.B. (1987 b).  Towards a unified asymptotic theory for
autoregression. {\em Biometrika 74}, 535-547.\\
Phillips P.C.B. (2007). Unit root log periodogram regression. {\em
Journal of Econometrics 138}, 104-124.\\
Resnick S. (1986). Point processes, regular variation and weak
convergence. {\em Advanced in Applied Probability 18}, 66-183.\\
Resnick S. (2007). {\em Heavy-Tail Phenomena.} Springer.\\
Rva\v{c}eva E.L. (1962). On domains of attraction of
multi-dimensional distributions. {\em Select. Transl. Math. Statist.
and Probability, Vol.2, American Mathematical Society, Providence,
RI}, 183-205.\\
Tyran-Kkami\'{n}ska M. (2010 a). Convergence to L\'{e}vy stable
processes under some weak dependence conditions. {\em Stochastic
Processes and their Applications 120}, 1629-1650.\\
Tyran-Kkami\'{n}ska M. (2010 b). Weak convergence to L\'{e}vy stable
processes in dynamical systems. {\em Stochastics and Dynamics  10},
263-289.

   \end{document}